\documentclass[11pt,psfig]{article}
\usepackage{amssymb,amsmath,amscd}
\usepackage[all]{xy}
\usepackage{mathtools}
\numberwithin{equation}{section}

\newtheorem{lem}{Lemma}[section]
\newtheorem{thm}{Theorem}[section]
\newtheorem{cor}{Corollary}[section]
\newtheorem{prop}{Proposition}[section]

\newcommand{\f}[1]{\mathfrak{#1}}

\newcommand{\commentout}[1]{}

\newcommand{\mc}{\mathcal}
\newcommand{\arr}[1]{\left( \begin{array}{clcr} #1 \end{array} \right)}

\newcommand{\sgn}{\, {\rm sgn}}

\DeclarePairedDelimiter\floor{\lfloor}{\rfloor}
\newcommand{\vertiii}[1]{{\left\vert\kern-0.25ex\left\vert\kern-0.25ex\left\vert #1 
    \right\vert\kern-0.25ex\right\vert\kern-0.25ex\right\vert}}

\begin{document}
\title{Certain $L^2$-norms on Automorphic Representations  of $SL(2, \mathbb R)$}
\author{Hongyu He  \\
Department of Mathematics \\
Louisiana State University \\
email: hhe@lsu.edu\\
}
\date{}
\maketitle

\abstract{Let $\Gamma$ be a non-uniform lattice in $SL(2, \mathbb R)$. In this paper, we study various $L^2$-norms of automorphic representations of $SL(2, \mathbb R)$. We bound these norms with intrinsic norms defined on the representation. Comparison of these norms will help us understand the growth of $L$-functions in a systematic way (~\cite{hers}). }

\section{Introduction}
Let $\Gamma$ be a non-uniform lattice in $SL(2, \mathbb R)$. By an automorphic representation of $SL(2, \mathbb R)$, we mean a   finitely generated admissible representation of $SL(2, \mathbb R)$, consisting of $\Gamma$-invariant functions on $SL(2, \mathbb R)$ (\cite{hc}). Among  all automorphic representations, $L^2$ automorphic representations, i.e., subrepresentations of $L^2(G/\Gamma)$, are of fundamental importance. Since $L^2$ automorphic representations are unitary and completely reducible, we assume $L^2$ automorphic representations to be irreducible.  By Langlands theory, $L^2$ automorphic representations come from either the residues of Eisenstein series or the cuspidal automorphic representations. Throughout this paper, we shall mostly focus on  irreducible  cuspidal representations, even though our results also apply to unitary Eisenstein series with vanishing constant term near a cusp.
 \\
\\
\commentout{ Let $f \in L_{loc}(G/\Gamma)$ satisfying the following conditions:
\begin{enumerate}
\item $f$ is cuspidal at $0$ and $\infty$, namely,
$$\int_{N_0/\Gamma_{0}} f(x n) d n=0, \qquad \int_{N_{\infty}/\Gamma_{\infty}} f (x n) d n=0.$$
Or equivalently, the zero-th Fourier coefficients with respect to $N_0$ and $N$ are zero.
\item $f(x w)=c_w f (x)$ where $c_w$ is a constant only depends on $w$.
\end{enumerate}
Denote the space of such functions by $L^0_{loc}(G/\Gamma)$.\\
\\
If  $f \in L^2(G/\Gamma)$ and 
 the  representation of $G$, generated by $f$, is an admissible quasisimple representation $\pi$ of $G$,
we write $f \in L^2(G/\Gamma)_{\pi}$. Such an automorphic representation will be called automorphic representation of type $\pi$.}

Let $G=SL(2, \mathbb R)$ and  $\pi$ be an irreducible admissible representation of $G$.   We say an automorphic representation is of type $\pi$ if the automorphic representation is infinitesimally equivalent to $\pi$. In particular, we write $L^2(G/\Gamma)_{\pi}$ for the sum of all $L^2$-automorphic representations of type $\pi$. It is well-known that $L^2(G/\Gamma)_{\pi}$ is of finite multiplicity (\cite{hc}). The main purpose of this paper is to study various $L^2$-norms of the automorphic forms at the representation level. In the literature, automorphic forms, the $K$-finite vectors in an automorphic representation, are the main focus of interests. Our main focus here is the $L^2$-norms of automorphic forms, in comparison with (intrinsic) norms in the representation. We hope to gain some understanding of various $L^2$-norms of  automorphic representation as a whole,  without references to automorphic forms. We believe this may lead to a better understanding of the Fourier coefficients and  $L$-functions. \\
\\
Our estimates of $L^2$-norms essentially involve two decompositions, the Iwasawa decomposition $KAN$, and its variant $KNA$. The $KAN$ decomposition is utilized mainly to define Fourier coefficients and constant terms of automorphic forms. We give estimates of various  $L^2$ norms of the restriction of automorphic representation to $AN$ and the Siegel set. The $KNA$ decomposition, on the other hand, seems to be a potentially useful tool to study the $L$-function associated with the automorphic representation. In this paper, we give various estimates on the $L^2$-norm of automorphic representation restricted to $\Omega A$, with $\Omega$ a compact domain  in $KN$. \\
\\
Our view point and setup are very similar to those of Harish-Chandra (\cite{hc}). The group action will be from the left and the standard cusp will be at zero instead of $\infty$. Working in the general framework of harmonic analysis on semisimple Lie groups,  Harish-Chandra gave a very detailed account of the theory of cusp forms and Eisenstein series, mainly due to Selberg, Gelfand and Piatetsky-Shapiro, and Langlands.  Our goal here is quite limited: we  only treat the group $G=SL(2, \mathbb R)$ and  we  study various $L^2$-norms of automorphic representations of type $\pi$. Most of our results are stated in terms of automorphic distribution (\cite{br} \cite{sc} \cite{ms1}). The reason is simple. There are two types of norms involved, one for the automorphic forms, and one for the representation. Using automorphic distributions, automorphic forms can be viewed as matrix coefficients of $K$-finite vectors and a fixed automorphic distribution. This allows us to compare norms of automorphic forms and norms of the representation. These results will shed lights on the growth of the Rankin-Selberg $L$-functions (\cite{hers}). \\
\\
To state our results in a simpler form, let $\Gamma=SL(2, \mathbb Z)$. Fix the usual Iwasawa decomposition $G=KAN$ with $N$ the unipotent upper triangular matrices. Let $\mc F$ be the fundamental domain of $G/\Gamma$ contained in  a Siegel set. Recall that the $L^2$-norm on the fundamental domain is
$$\| f \|_{L^2(G/\Gamma)}^2= \int_{\mc F} |f(kan)|^2 a^{2} \frac{ d a}{a} d n dk.$$
We have
\begin{thm}Let $\pi=\mc P(u, \pm)$ be a unitary representation in the principal series (see Section \ref{psr} for the definition).
Let $\mc H$ be a cuspidal representation in $L^2(G/\Gamma)_{\pi}$. Then  for any  $\epsilon >0$, there exists a $C_{\epsilon} >0$ such that
$$\int_{{\mc F}} | f (kan)|^2 a^{\epsilon} \frac{d a}{a} d n d k  \leq C_{\epsilon} 
\| f\|_{L^2(G/\Gamma)}^2, (\forall \,\,\, f \in \mc H).$$
For any $\epsilon <0$, there exists a $C_{\epsilon}>0$ such that
$$\int_{{\mc F}} | f (kan)|^2 a^{\epsilon} \frac{d a}{a} d n d k  \leq C_{\epsilon} 
\vertiii{f}_{\frac{\epsilon}{2}-u_0}, (\forall \,\,\, f \in \mc H^{\infty}).$$
Here $u_0=\Re(u)$ and the norm $\vertiii{f}_{\frac{\epsilon}{2}-u_0}$ is defined on ${\mc H}^{\infty}$, smooth vectors in the representation in $\mc H$( see Eq. \ref{triplenorm} for the definition of $\vertiii{f}$).
\end{thm}
Our theorem essentially says that every $f \in L^2(G/\Gamma)_{\pi}$ is also in 
$L^2(\mc F, a^{\epsilon} \frac{  d a}{a} d n d k)$ for every $\epsilon>0$. In other words, the natural injection
$$L^2(\mc F, a^2 \frac{  d a}{a} dn d k) \supseteq \mc H \rightarrow L^2(\mc F, a^{\epsilon} \frac{  d a}{a} dn d k)$$
is bounded for every $\epsilon >0$ even though the natural map 
$$L^2(\mc F, a^2 \frac{  d a}{a} dn d k)  \rightarrow L^2(\mc F, a^{\epsilon} \frac{  d a}{a} dn d k)$$
is not bounded unless $\epsilon \geq 2$. In terms of the parameter $\epsilon$,
there is a natural barrier at $\epsilon=0$, namely, as $\epsilon \rightarrow 0$, the norms of these bounded operators go to infinity.\\
\\
We shall remark that our estimates are true for all nonuniform lattices of any finite covering of $SL(2, \mathbb R)$ (see Theorem \ref{main1}). In addition, the first bound with $\epsilon >0$ also holds for discrete series $D_n$ (see Cor. \ref{discrete}). They are proved by studying the $L^2$-norms of Fourier coefficients of the automorphic distribution, defined in Schmid (\cite{sc}) and Bernstein-Resnikov (\cite{br}).  For the general linear group $GL(n, \mathbb R)$, similar results should hold. The following  problem is worthy of further investigation.\\
\\
{\bf Problem}: Let $G$ be a semisimple Lie group,  $\Gamma$ an arithmetic lattice and $S$ a Siegel domain. Find the best exponents $\alpha$ such that
$$i: L^2(G/\Gamma)_{\pi} \rightarrow L^2(S, a^{\alpha} \frac{d a}{a} d n dk)$$
is bounded. Here $G=KAN$ is the Iwasawa decomposition.\\
\\
Notice that if $\alpha= 2 \rho$, the sum of positive roots of $\f{ gl}(n)$, the measure on the right hand side is the invariant measure of $G$ restricted to $S$. In this case, $i$ is automatically bounded. This shows that if $\alpha$ is \lq\lq bigger" than $2 \rho$, $i$ is also bounded. The problem is to find the \lq\lq smallest" $\alpha$ such that $i$ is bounded. We shall remark that cusp forms will remain to be in $L^2(S, a^{\alpha} \frac{d a}{a} d n dk)$ for any $\alpha$ since they are fast decaying on the Siegel set. Hence our problem is about cuspidal representations, rather than cusp forms.\\
\\
The second main result is an $L^2$-estimates of $f$ on $\Omega A$ where $\Omega$ is a compact domain in $G/A$.
\begin{thm}
Let $\Gamma$ be a nonuniform lattice  in $SL(2, \mathbb R)$. Suppose that the Weyl element $w \in \Gamma$ and $\Gamma \cap N \neq \{ I \}$. Let $\mc H$ be a cuspidal  automorphic representation of $G$ of type $\mc P(i \lambda, \pm)$.   Let $\Omega$ be a compact domain in $KN$.  Let $\epsilon >0$. Then there exists a positive constant $C$ depending on $\epsilon, \mc H$ and $\Omega$ such that
$$\| f \|_{L^2(\Omega A, a^{\epsilon} \frac{ d \, a }{a} d t d k)} \leq C \vertiii{f}_{-\frac{\epsilon}{2}} \qquad (f \in \mc H^{\infty}).$$
See Eq. \ref{triplenorm} for the definition of $\vertiii{f}$.
\end{thm}
We shall remark that in the $KNA$ decomposition, the invariant measure is given by $d k d n \frac{da}{a}$. Hence, the $L^2$-norm here is a perturbation of the canonical $L^2$-norm. In addition,  $\Omega A$ has infinite measure. The perturbation is needed because our theorem fails at $\epsilon=0$. At $\epsilon=0$, the norm  $\vertiii{f}_{-\frac{\epsilon}{2}}$ is the original Hilbert norm  $\| f\|$ of the cuspidal representation. There is no  chance that
$\| f \|_{L^2(\Omega A, \frac{ d \, a }{a} d t d k)} $ can remain bounded for all $f \in \mc H$. \\
\\
Throughout our paper, the Haar measure on $A$ will be $\frac{ d a}{a}$. We use $c$ or $C$ as  symbolic constants and $c_{\epsilon, u}$ to indicate the dependence on $\epsilon$ and $u$.

\section{Certain $L^2$-norm of $\Gamma$-invariant functions}
\subsection{setup}
Let $G=SL(2, \mathbb R)$.
Let $$N=\{ n_t=\arr{ 1 & t  \\ 0 & 1} : t \in \mathbb R \},$$
 $$K=\{ k_{\theta}=\arr{\cos \theta & -\sin \theta \\ \sin \theta & \cos \theta}: \theta \in [0, 2\pi) \}$$
 $$A=\{ \arr{a & 0 \\ 0 & a^{-1}} : a \in \mathbb R^+\},$$
 and $w =\arr{0 & 1 \\ -1 & 0} \in K$. We call $w$ the Weyl element. Let $\Gamma$ be a discrete subgroup of $G$ such that 
$\Gamma \cap N $ is nontrivial. Without loss of generality assume that 
$$\Gamma \cap N=N_{\mathfrak p}=\{ n_t: t \in \mathfrak p \mathbb Z \}$$
with $\mathfrak p \in \mathbb N^+$.\\
\\

Let $M= \{ \pm I \} \subseteq K$. Fix $P=MAN$, the minimal parabolic subgroup. Then the identity component  $P_0=AN$. Fix $\frac{ d a}{a} d t$ as the left invariant measure on $an_t \in P_0$ and $d T \frac{ da }{a} $ as the right invariant measure on $N_{T} a \in P_0$. We shall keep the notion that  $a n_t= N_{T} a$. Then
$$ T=a^2 t, \qquad t=a^{-2} T, \qquad   \frac{ d a}{a} d T= a^2 \frac{ d a}{a} d t, \qquad \frac{ d a}{a} d t=a^{-2} \frac{ d a}{a} d T.$$
Fix $ d k= d \theta $ as the invariant measure on $K$. We write $g=k_{\theta} a n_t$ for the $KAN$ decomposition and $g=k_{\theta} n_T a $ for  the $KNA$ decomposition. Fix the standard invariant measure
$$ d g= a^2  d t \frac{d a}{ a} d k=  d T \frac{ d a}{a} d k.$$\\
\\ 
Let $N_{T_1}=\{ n_T: 0 \leq T \leq T_1 \}$ if $T_1>0$ and $N_{T_1}=\{ n_T: 0 \geq T \geq T_1 \}$ if $T_1 <0$. Let $X_{T_1}=K N_{T_1} A$ equipped with the canonical measure $d k d T \frac{ d a}{a}$. Let $\epsilon \in \mathbb R$. For $f \in C(G/\Gamma)$ or more generally $L^2_{loc}(G/\Gamma)$, we would like to estimate
$$\|f \|_{T_1, \epsilon}=\| f \|_{L^2(X_{T_1}, a^{\epsilon} \frac{d a}{a} d k  d T)}.$$
Here $L^2_{loc}(G/\Gamma)$ is the space of locally square integrable function on $G/\Gamma$.
\\
\\
Let $a_1 \in \mathbb R^+$. Let $A_{a_1}^+=\{ a \geq a_1 \}$ and $A_{a_1}^-=\{ 0 < a \leq a_1\}$. By abusing notation, we simply use $a \in \mathbb R^+$ as an element in $A$. Write
$$X(T_1, a_1)^{\pm}= K N_{T_1} A_{a_1}^{\pm} , \qquad P(T_1, a_1)^{\pm}= N_{T_1} A_{a_1}^{\pm} .$$
Write $\| f \|_{L^2(X(T_1, a_1)^{\pm}, a^{\epsilon} \frac{d a}{a}d T d k)}$ as $ \| f \|_{T_1,  a_1^{\pm},\epsilon}$.
\subsection{Estimates on $\| f \|_{T_1,  a_1^{-},\epsilon}$}
Without loss of generality, assume
$T_1 >0$. 
Observe  that $$P(T_1, a_1)^{-}= \{ 0 \leq T \leq T_1, 0< a \leq a_1\}= \{ 0 < a \leq a_1, 0 \leq t \leq a^{-2} T_1\}.$$
We have
\begin{prop} Let $f \in L^2_{loc}(P_0)$ such that $f(x N_{\mathfrak p})=f(x)$ for a fixed period $\mathfrak p \in \mathbb N^+$. Then for any $\epsilon \in \mathbb R$,
$$  \int_0^{a_1} a^{2+\epsilon} \floor{ \frac{T_1}{{\mathfrak p} a^2}} \int_0^{\mathfrak p} | f ( a n_t)|^2 d t \frac{ d a}{a} \leq \|f\|_{L^2(P(T_1,a_1)^-, a^{\epsilon} \frac{ d a}{a} d T)}^2 \leq 
\int_0^{a_1} a^{2+\epsilon} (\floor{ \frac{T_1}{{\mathfrak p} a^2}}+1) \int_0^{\mathfrak p} | f (a n_t)|^2 d t \frac{ d a}{a}.$$
\end{prop}
Proof: We have
\begin{equation}
\begin{split}
 & \|f\|^2_{L^2(P(T_1,  a_1)^{-}, a^{\epsilon} \frac{d a}{a} d T)} \\
 = & \int_{0}^{a_1} \int_0^{T_1} a^{\epsilon} \| f( n_T a) \|^2  d T  \frac{ d a}{ a}\\
 = & \int_{0}^{a_1} \int_0^{a^{-2} T_1} a^{2+\epsilon} \| f(a n_t) \|^2 dt \frac{d a}{ a} \\
 \geq & \int_{0}^{a_1} a^{2+\epsilon} \floor{\frac{T_1}{a^2 {\mathfrak p}}} (\int_0^{{\mathfrak p}}  \| f(a n_t) \|^2 dt )\frac{d a}{ a} 
\end{split}
\end{equation}
Here $\floor{*}$ is the floor function. The other direction is similar. $\Box$ \\
\\
For $T_1$ negative, we have a similar statement. Combining these two cases, we have
\begin{thm}\label{thm1.1}
Assume that $f \in L_{loc}^2(G)$ and $f(x N_{\mathfrak p})=f(x)$ for a fixed period ${\mathfrak p}$. Let $a_1 > 0$ and $\epsilon \in \mathbb R$. Then
$$ \int_K \int_0^{a_1} a^{2+\epsilon} \floor{ \frac{|T_1|}{{\mathfrak p} a^2}} \int_0^{\mathfrak p} | f (k a n_t)|^2 d t \frac{ d a}{a} d k \leq \|f\|_{T_1, a_1^{-},\epsilon }^2 \leq 
 \int_K \int_0^{a_1} a^{2+\epsilon} (\floor{ \frac{|T_1|}{{\mathfrak p} a^2}}+1) \int_0^{\mathfrak p} | f(k a n_t)|^2 d t \frac{ d a}{a} d k.$$
\end{thm}

\subsection{Estimate on  $\| f \|_{T_1,  a_1^{+}, \epsilon}$}
To estimate $\| f \|_{T_1, a_1^{+},\epsilon}$, we must utilize the Weyl group element $w$. We  assume that 
$$| f( x w) |= | f(x)| \qquad  (\forall \, x \in G). $$ 
Let $a \in [ a_1, \infty)$. By the Iwasawa decomposition
$$n_T a w=k(T,a) n_{T^{\prime}} a^{\prime}, \qquad a^{\prime}=\frac{\sqrt{T^2+1}}{a}, \,\, T^{\prime}=-T$$
and $k(T,a) \in K$. This defines a coordinate transform from $(T, a)$ to $(T^{\prime}, a^{\prime})$. 
Let $(P({T_1, a_1})^+)^{\prime}$ be the coordinate transform of $P(T_1, a_1)^+ w$ in terms of $(T^{\prime}, a^{\prime})$ coordinates. We have
$$(P(T_1, a_1)^+)^{\prime}= \{ -T_1 \leq T^{\prime} \leq 0, 0 < a^{\prime} \leq \frac{\sqrt{(T^{\prime})^2+1}}{a_1} \}.$$
It is easy to see that
$$P(-T_1, \frac{1}{a_1})^- \subseteq (P(T_1, a_1)^+)^{\prime} \subseteq P(-T_1, \frac{\sqrt{T_1^2+1}}{a_1})^-,$$
and
$$K P(-T_1, \frac{1}{a_1})^- \subseteq K P(T_1, a_1)^+ w \subseteq K P(-T_1, \frac{\sqrt{T_1^2+1}}{a_1})^-.$$
Observe that $$a^{\epsilon} \frac{d a}{a} d T=(\sqrt{(T^{\prime})^2+1})^{\epsilon} (a^{\prime})^{-\epsilon} d T^{\prime} \frac{d a^{\prime}}{a^{\prime}}$$
 and $$|f(k n_T a)|^2= |f(k n_T a w)|^2= |f(k k(T,a) n_{T^{\prime}} a^{\prime})|^2.$$  We obtain
\begin{prop} Let $f \in L^2_{loc}(G)$, $a_1 >0$ and $\epsilon \in \mathbb R$. Suppose that $f(x N_{\mathfrak p})=f(x)$ and $| f(xw)|= |f(x)|$. Then 
$$ \| f \|^2_{-T_1, (\frac{1}{a_1})^{-}, -\epsilon} \leq \|f \|_{T_1, a_1^{+}, \epsilon}^2 \leq 
(\sqrt{T_1^2+1})^{\epsilon} \| f\|^2_{-T_1, (\frac{\sqrt{T_1^2+1}}{a_1})^{-},  -\epsilon} \qquad (\epsilon \geq 0);$$
$$ (\sqrt{T_1^2+1})^{\epsilon} \| f \|^2_{-T_1, (\frac{1}{a_1})^{-}, -\epsilon} \leq \|f \|_{T_1, a_1^{+}, \epsilon}^2 \leq 
\| f\|^2_{-T_1, (\frac{\sqrt{T_1^2+1}}{a_1})^{-},  -\epsilon} \qquad (\epsilon \leq 0);$$
\end{prop}

\subsection{Estimates of $\|f\|_{T_1, \epsilon}$}
Choose $a_1=1$. We have 
$$ \| f \|^2_{-T_1, (1)^{-}, -\epsilon} \leq \|f \|_{T_1, 1^{+}, \epsilon}^2 \leq 
(\sqrt{T_1^2+1})^{\epsilon} \| f\|^2_{-T_1, (\sqrt{T_1^2+1})^{-},  -\epsilon} \qquad ( \epsilon \geq 0);$$
$$ (\sqrt{T_1^2+1})^{\epsilon} \| f \|^2_{-T_1, (1)^{-}, -\epsilon} \leq \|f \|_{T_1, 1^{+}, \epsilon}^2 \leq 
 \| f\|^2_{-T_1, (\sqrt{T_1^2+1})^{-},  -\epsilon} \qquad ( \epsilon \leq 0);$$
Combined with Theorem \ref{thm1.1}, we have
\begin{thm} Let $f$ be a locally square integrable function on $SL(2, \mathbb R)$ such that $f(x N_{\mathfrak p})=f(x)$ and $| f(xw)|= |f(x)|$. If $\epsilon >0$,  then
\begin{equation}
\begin{split}
& \int_K \int_0^1 (a^{2+\epsilon}+a^{2-\epsilon}) \floor{\frac{T_1}{{\mathfrak p}a^2}} \int_0^{\mathfrak p} |f(k a n_t)|^2 d t \frac{ d a}{a} d k \leq  \| f \|_{T_1, \epsilon}^2 \\
\leq 
  & \int_K \int_0^1 a^{2+\epsilon}  (\floor{\frac{T_1}{{\mathfrak p}a^2}}+1) \int_0^{\mathfrak p} |f(k a n_t)|^2 d t \frac{ d a}{a} d k 
  +  (\sqrt{T_1^2+1})^{\epsilon} \int_K  \int_{0}^{\sqrt{T_1^2+1}} a^{2-\epsilon}(\floor{\frac{T_1}{{\mathfrak p}a^2}}+1) \int_0^{\mathfrak p} |f(k a n_t)|^2 dt \frac{ da}{a} dk .
\end{split}
\end{equation}
If $\epsilon \leq 0$, then 
\begin{equation}
\begin{split}
& \int_K \int_0^1 (a^{2+\epsilon}+(\sqrt{T_1^2+1})^{\epsilon} a^{2-\epsilon}) \floor{\frac{T_1}{{\mathfrak p}a^2}} \int_0^{\mathfrak p} |f(k a n_t)|^2 d t \frac{ d a}{a} d k \leq  \| f \|_{T_1, \epsilon}^2 \\
\leq 
  & \int_K \int_0^1 a^{2+\epsilon}  (\floor{\frac{T_1}{{\mathfrak p}a^2}}+1) \int_0^{\mathfrak p} |f(k a n_t)|^2 d t \frac{ d a}{a} d k 
  +  \int_K  \int_{0}^{\sqrt{T_1^2+1}} a^{2-\epsilon}(\floor{\frac{T_1}{{\mathfrak p}a^2}}+1) \int_0^{\mathfrak p} |f(k a n_t)|^2 dt \frac{ da}{a} dk .
\end{split}
\end{equation}
\end{thm}
If $\mathfrak p=1$ and $T_1=1$, we have
$$ \| f \|_{T_1, \epsilon}^2  \leq C_{\epsilon} \int_K \int_0^{\sqrt{2}} (a^{\epsilon}+a^{-\epsilon}) \int_0^1 | f (ka n_t)|^2 d t \frac{ d a}{ a} d k.$$
Notice that  for $ 0 < a \leq \sqrt{2}$, $\floor{\frac{1}{a^2}}+ 1 \leq \frac{2}{a^2}$.
Hence, we have bounded the norm of $f$ on $X_{T_1}$. Generally, we have

\begin{thm}\label{main}Suppose that $f$ is a locally square integrable function on $SL(2, \mathbb R)$ such that $f(x N_{\mathfrak p})=f(x)$ and $| f(xw)|= |f(x)|$. Let $\epsilon \in \mathbb R $. Then there exists a positive constant $c_{T_1,\epsilon, \mathfrak p}$ such that
$$ \| f \|_{T_1, \epsilon}^2  \leq c_{T_1,\epsilon, \mathfrak p} \int_K \int_0^{\sqrt{1+T_1^2}} ( a^{\epsilon}+a^{-\epsilon}) \int_0^{\mathfrak p} | f (ka n_t)|^2 d t \frac{ d a}{ a} d k.$$
\end{thm}
Proof: We choose a positive constant $c$ such that
$$\floor{\frac{T_1}{\mathfrak p a^2}}+1 \leq c \frac{T_1}{\mathfrak p a^2} \qquad (\,\, \forall \,\,\,\, 0 < a \leq \sqrt{1+T_1^2}).$$
Then let $c_{T_1,\epsilon, \mathfrak p}=c \max(2, 1+(\sqrt{1+T_1^2})^{\epsilon}) $. 
$\Box$ \\
\\
Observe that the right hand side of our inequality involves an integral over a Siegel set. However the measure on this Siegel set can be larger than the invariant measure $ a^2 d k \frac{ d a}{a} d t$. What we have achieved is a bound of $\| f \|_{T_1, \epsilon}$ by an integral on a Siegel set. In the next section, we shall give estimation of the norms of $f$ on $A_{a_1}^{-} N/N_{\mathfrak p}$ and on $KA_{a_1}^{-} N/N_{\mathfrak p}$.

\section{Matrix Coefficients and Analysis on $P_0/N_{\f p}$}
Now we shall focus on $L^2$ automorphic representations of type $\pi$ where $\pi$ is a principal series representation. According to Langlands, $L^2$ automorphic representations come from either the residue of Eisenstein series or cuspidal automorphic forms. In either cases, the restrictions of $L^2$ automorphic representations fail to be $L^2$ on $P_0/N_{\f p}$, when $P_0/N_{\f p}$ is equipped with the left invariant measure. However if we perturb the invariant measure correctly, automorphic forms will be square integrable. In this section, we will discuss the $L^2$-integrability of $f|_{P_0}$ with $f \in L^2(G/\Gamma)_{\pi}$ with respect to the measure $a^{\epsilon} \frac{ d a} {a} d t$. We will consequently discuss the $L^2$-norm  on a Siegel subset. We conduct our discussion in terms of matrix coefficients with respect to periodical distributions with no constant term. More precisely, the function $f|_{P_0}$ will be regarded as the matrix coefficient of $v \in \mc H_{\pi}$ and a periodical distribution in $(\mc H^{*})^{-\infty}$. Our view is similar to Schmid and Bernstein-Reznikov (\cite{sc} \cite{br}).

\subsection{Principal series representations of $SL(2, \mathbb R)$}\label{psr}
Principal series representations of $G$ can be easily constructed using homogeneous distributions on $\mathbb R^2-\{0\}$, namely, those
$$\{ f(rx)=r^{-1-u} f(x), f(-x)=\pm f(x) \mid r \in \mathbb R^+, f \in C(\mathbb R^2 -\{0\}) \}. $$
See for example \cite{cas} \cite{knapp}. In this section, we shall focus on the smooth vectors and the space of distributions associated with them. Let $(\pi_{u, \pm}, \mathcal P(u, \pm))$ be the {\bf unitarized} principal series representation with the trivial or nontrivial central character. $\mc P(u, \pm )$ includes unitary principal series $\mc P(u, \pm)$ (with $u \in i \mathbb R$) and complementary series $\mc P(u,+)$ (with $u \in (-1,0) \cup (0,1)$). All of these representations are irreducible except $\mc P(0, -)$. In addition $\mc P(u, \pm) \cong \mc P(-u, \pm)$. \\
\\
Consider the noncompact picture (\cite{knapp}). The noncompact picture is essentially the restriction of $f$ onto the line $\{(x,1) \mid x \in \mathbb R \} \subseteq \mathbb R^2$. We have for any $g=\arr{a & b \\ c & d}$, $f \in \mathcal P(u, \pm)^{\infty}$,
$$\pi_{u, \pm}(g) f (x)=\chi_{\pm}(a-cx) |a - c x|^{-1-u} f(\frac{d x-b}{a-cx}).$$
Here $\chi_{-}(x)$ is the sign character on $\mathbb R - \{0 \}$ and $\chi_{+}(x)$ is the trivial character.
In particular, we have
$$\pi_{u, \pm}\arr{a & 0 \\ 0 & a^{-1}} f(x)=|a|^{-1-u} f(a^{-2} x), \qquad (a \in \mathbb R^+);$$
$$\pi_{u, \pm} \arr{1 & b \\ 0 & 1} f(x)= f(x-b);$$
$$\pi_{u, \pm} (w) f(x)=\chi_{\pm}(-x) |x|^{-1-u} f(-\frac{1}{x});$$
$$\pi_{u, \pm} \arr{\cos \theta & - \sin \theta \\ \sin \theta & \cos \theta} f(x)=\chi_{\pm}(\cos \theta- x \sin \theta) |\cos \theta- x \sin \theta|^{-1-u} f(\frac{ x \cos \theta+ \sin \theta}{\cos \theta- x \sin \theta}).$$
There is a $G$-invariant pairing between $\mc P(u, \pm)^{\infty}$ and $\mc P(-u, \pm)^{\infty}$. This allows us to write the dual space of $\mc P(u, \pm)^{\infty}$ as $\mc P(-u, \pm)^{-\infty}$.\\
\\
Unless otherwise stated, $\mc P(u, \pm)$ will refer to the noncompact picture.  The space $\mc P(u, \pm)^{\infty}$ will then be a subspace of infinitely differentiable functions on $N \cong \mathbb R$ satisfying certain conditions at infinity.
\subsection{ Matrix coefficients with respect to periodical distribution with zero constant term}
According to \cite{br} \cite{sc} \cite{ms}, every $L^2$ automorphic form of type $\pi$ can be written as matrix coefficients of an automorphic distribution and a  vector in the unitary representation $\pi$. Equivalently, in our setting,
there exists a distribution $\tau \in \mc P(u, \pm)^{-\infty}$ such that the automorphic forms of type $\pi$ can be written as  linear combinations of
$$f_m(g)=\langle \pi_{u, \pm}(g) \tau, v_m \rangle,$$
with $v_m(x)=(1+x^2)^{-\frac{1-u}{2}}(\frac{1+x i}{1-xi})^{\frac{m}{2}}$. For $\mc P(u, +)$, the weight $m$ can only be an even integer. For $\mc P(u, -)$, the weight $m$ must be an odd integer. If $\tau$ is cuspidal,  $\tau$ has a Fourier expansion
$$\tau=\sum_{n \in {\mathfrak p}^{-1} \mathbb Z, n \neq 0}^* b_{n} \exp  2 \pi i x n,$$
Here $\mathfrak p$ is a positive integer and $\sum^*$ denote the weak summation (\cite{he}). We call such  $\tau$  a periodical distribution without constant term. \\
\\
Let $\tau \in \mc P(u, \pm)^{-\infty}$ be a periodic distribution without constant term. We compute the matrix coefficient formally:
\begin{equation}\label{f(an)}
\begin{split}
&  \langle \pi_{u, \pm}(a n_t) \tau, v \rangle \\
 = & \langle \sum_{n \in {\mathfrak p}^{-1} \mathbb Z, n \neq 0} a^{-1-u} b_{n} \exp 2 \pi i (a^{-2}x-t) n, v(x) \rangle \\
 = & a^{-1-u} \sum_{n \in {\mathfrak p}^{-1} \mathbb Z, n \neq 0}^* \int b_{n} \exp (2 \pi i a^{-2} x n) \exp (-2 \pi i t n ) v(x) d x \\
 = & a^{-1-u} \sum_{n \in {\mathfrak p}^{-1} \mathbb Z, n \neq 0}^*  b_{n} (\mc F v)(-n a^{-2} )  \exp (- 2 \pi i  t n).
 \end{split}
 \end{equation}
 Here  $\mc F$ is the Fourier transform, and $v$ is in a suitable subspace of $ \mc P(-u,\pm)^{-\infty}$.  The formula above, also known as the Fourier-Whittaker expansion in a more general context, is valid for
 $v \in \mc P(-u, \pm)^{\infty}$ with $\Re(-u)>-1$.
 
 \begin{lem}\label{nsum} Let $u=u_0+i u_1 $ with $u_0 <1$ and 
 $$\tau =\sum_{n \in {\mathfrak p}^{-1} \mathbb Z, n \neq 0}^* b_{n} \exp  2 \pi i x n \in \mc P(u, \pm)^{-\infty}.$$ 
 For $v \in \mc P(-u, \pm)^{\infty}$, we have
  $$\langle \pi_{u, \pm}(a n_t) \tau, v \rangle=a^{-1-u} \sum_{n \in {\mathfrak p}^{-1} \mathbb Z, n \neq 0}^*  (\mc F v)(-n a^{-2} )  b_{n} \exp (- 2 \pi i  t n)$$
 $$\int_0^{\mathfrak p} |\langle \pi_{u, \pm}(a n_t) \tau, v \rangle|^2 d t= {\mathfrak p} \sum_{{n \in {\mathfrak p}^{-1} \mathbb Z}} a^{-2-2u_0} |b_{n}|^2 | \mc F v(-n a^{-2})|^2.$$
\end{lem}
Proof:  Suppose $\Re(u) <1$. The functions in $\mc P(-u, \pm)^{\infty}$ are smooth functions of the form $(1+x^2)^{-\frac{1-u}{2}} \phi(\frac{1+xi}{1-xi})$ with $\phi$ an odd or even smooth function on the unit circle. They are slowly decreasing functions. Their Fourier transforms exist. Since the derivatives $v^{(n)}$  are  of this form and they are integrable , we see that $\mc F v(\xi)$ will decay faster than any polynomial at $\infty$. The weak sum in Equation (\ref{f(an)}) becomes a convergent sum. Our lemma is proved. $\Box$\\
\\
We shall make a few remarks here. Since $v \in \mc P(-u, \pm)^{\infty}$ and $\tau \in \mc P(u, \pm)^{-\infty}$, the matrix coefficient  $\langle \pi_{u, \pm}(a n_t) \tau, v \rangle$ is automatically smooth. Our lemma simply provided a Fourier expansion, which is generally known as the Fourier-Whittaker expansion over the whole group $G$. 
The restriction that $u_0<1$ is somewhat  unsatisfactory. When $u_0 \geq 1$, $\mc Fv(\xi)$ may fail to be a function even for $v$ smooth. This happens when $\mathcal P(-u, \pm)$ is reducible and discrete series will appear as composition factors. Hence, automorphic representations that are  discrete series, can be treated by considering the reducible $\mc P(-u, \pm)$. We shall refer readers to Schmid's paper \cite{sc} for details.  When $\mathcal P(-u, \pm)$ is irreducible, $\mc Fv(\xi)$ is a fast decaying continuous function off from zero. Our lemma is still valid in this case. However,if $u_0 >1$,  $\mc Fv(\xi)$ will fail to be a locally integrable function near zero and need to be regularized to be a Schwartz distribution.  \\
\\
 From now on, without further mentioning,  we will restrict our scope to $u_0<1$. We do not lose any generalities here. If $\mc P(u, \pm)$ is unitary, then $\Re(u) \in (-1, 1)$. If $\pi$ is a discrete series representation, then $\pi$ can be embedded into a principal series representation $\mc P(-u, \pm)$ with $u <1$. Hence  our assumption is adequate for the discussion of  $L^2$ automorphic representations. When $\Re(u) < 1$ and $v \in \mc P(-u, \pm)^{\infty}$,  $\langle \exp 2 \pi i x n, v \rangle$ shall be interpreted as 
$$-\frac{1}{ 2 \pi i n} \langle \exp 2 \pi i n x, \frac{d v}{ d x} \rangle.$$

 \subsection{$L^2$-norms on $P_0/N_{\f p}$}

Let us first study the $L^2$ norms of $f(g)=\langle \pi_{u, \pm}(g) \tau, v \rangle$ on $P_0/N_{\f p}$. $\tau$ and $v$ are given in Lemma \ref{nsum}.
Now we compute
 \begin{equation}\label{fc}
 \begin{split}
  & \int_{0}^{a_1} \int_0^{\mathfrak p} |f(a n_t)|^2 d t a^{\epsilon} \frac{ d a}{a} \\
=  & {\mathfrak p} \int_0^{a_1} a^{\epsilon} \sum_{n \in {\mathfrak p}^{-1} \mathbb Z, n \neq 0} a^{-2-2u_0} |b_{n}|^2 | \mc Fv(-n a^{-2})|^2 \frac{d a}{a} \\
= & {\mathfrak p} \int_{a_1^{-1}}^{\infty}  a^{-\epsilon} \sum_{{n \in {\mathfrak p}^{-1} \mathbb Z},n \neq 0} a^{2+2u_0} |b_{n}|^2 | \mc Fv(-n a^{2})|^2 \frac{d a}{a} \\
= & \frac{{\mathfrak p}}{2} \sum_{n \in {\mathfrak p}^{-1} \mathbb Z, n \neq 0}  \int_{a_1^{-2}}^{\infty}  a^{-\frac{\epsilon}{2}+1+u_0}   |b_{n}|^2 | \mc Fv(-n a)|^2 \frac{d a}{a} \\
= & \frac{{\mathfrak p}}{2} \sum_{{n \in {\mathfrak p}^{-1} \mathbb Z}, n > 0} \sum_{ \pm } \int_{\frac{n}{a_1^2}}^{\infty} a^{-\frac{\epsilon}{2}+1+u_0} n^{\frac{\epsilon}{2}-1-u_0} |b_{ \pm n}|^2 | \mc Fv( \mp a)|^2 \frac{d a}{a} \\
= & \frac{\mathfrak p}{2} \sum_{\pm} \int_{\frac{1}{a_1^2 \mathfrak p}}^{\infty} a^{-\frac{\epsilon}{2}+u_0} | \mc Fv( \mp a)|^2 \left[ \sum_{\frac{1}{\mathfrak p} \leq n \leq a  a_1^2, {n \in {\mathfrak p}^{-1} \mathbb Z}} n^
{\frac{\epsilon}{2}-1-u_0} |b_{\pm n}|^2 \right] d a \\
 \end{split}
 \end{equation}
 We summarize this in the following proposition.
 
 \begin{prop}\label{pestimate} Let $u=u_0+i u_1$ with $u_0 <1$. Let $v \in \mathcal P(-u, \pm)^{\infty}$ and $\tau \in \mathcal P(u, \pm)^{-\infty}$:
 $$\tau=\sum_{n \in {\mathfrak p}^{-1} \mathbb Z, n \neq 0}^* b_n \exp (2 \pi i n x).$$
 Let $f(an_t)=\langle \pi_{u, \pm}(a n_t) \tau, v \rangle$. Then $f(a n_t)$ is a smooth function on $P_0$ and 
 \begin{equation}
 \int_{0}^{a_1} \int_0^{\mathfrak p} |f(a n_t)|^2 d t a^{\epsilon} \frac{ d a}{a}=\frac{\mathfrak p}{2} \sum_{\pm} \int_{\frac{1}{a_1^2 \mathfrak p}}^{\infty} a^{-\frac{\epsilon}{2}+u_0} | \mc Fv( \mp a)|^2 \left[\sum_{\frac{1}{\mathfrak p} \leq n \leq a  a_1^2, {n \in {\mathfrak p}^{-1} \mathbb Z}} n^
{\frac{\epsilon}{2}-1-u_0} |b_{\pm n}|^2 \right] d a.]
\end{equation}
In particular,
\begin{equation}\label{p0es}
\int_{0}^{\infty} \int_0^{\mathfrak p} |f(a n_t)|^2 d t a^{\epsilon} \frac{ d a}{a}=\frac{\mathfrak p}{2} \sum_{\pm}  \left[\sum_{\frac{1}{\mathfrak p} \leq n, {n \in {\mathfrak p}^{-1} \mathbb Z}} n^
{\frac{\epsilon}{2}-1-u_0} |b_{\pm n}|^2 \right] \int_{0}^{\infty} a^{-\frac{\epsilon}{2}+u_0} | \mc Fv( \mp a)|^2 d a.
\end{equation} 
\end{prop}
Proof: Since $f(g)$ is a smooth function on $G$,  $f(a n_t)$ is a smooth function on $P_0$. Both equations hold without any assumptions on convergence.
Hence both sides of the equations converge or diverge at the same time. $\Box$

\subsection{Estimates of Fourier coefficients $b_{n}$}
We can now provide some estimates of certain sum of Fourier coefficients. These estimates are more or less known for automorphic forms (\cite{br} \cite{sc} \cite{ms1} \cite{go}). Our setting is more general.
 \begin{thm}\label{fcoef} Under the same assumption as Prop. \ref{pestimate}, suppose that there exists a $v \in \mc P(-u, \pm)^{\infty}$ such that $f(an_t)=\langle \pi_{u, \pm}(an_t) \tau, v \rangle$  is bounded on $P_0$. Suppose that $\mc F v(a)$ is nonvanishing on $\mathbb R^{-}$ or $\mathbb R^{+}$. Then we have the following estimates about the Fourier coefficients 
 $b_{n}$.
 \begin{enumerate}
 \item If 
 $|f(a n_t)|^2 \leq C_{\mu,f} a^{\mu }$ for some $\mu>0$, i. e., $f(an_t)$ decays faster than $a^{\mu}$ near the cusp $0$,  then  we have for each $\epsilon \in (-\mu, 0)$,
 $$
 \sum_{n>0,{n \in {\mathfrak p}^{-1} \mathbb Z} } n^{\frac{\epsilon}{2}-1-u_0} |b_{ \pm n}|^2 < \infty.$$
 \item For each $\epsilon>0$, there exists a $C_{\epsilon, \tau}>0$ such that 
 $$\sum_{n=\frac{1}{\mathfrak p}, {n \in {\mathfrak p}^{-1} \mathbb Z}}^k n^{\frac{\epsilon}{2}-1-u_0} |b_{ \pm n}|^2 < C_{\epsilon,\tau} k^{\frac{\epsilon}{2}} \qquad (k>1).$$
\commentout{ \item There exists a $C_{0, \tau} > 0$ such that
  $$\sum_{n=\frac{1}{\mathfrak p}, {n \in {\mathfrak p}^{-1} \mathbb Z}}^{k} n^{-1-u_0} |b_{\pm n}|^2 < C_{0,\tau} \ln k \qquad (k>1).$$}
 \end{enumerate}
 \end{thm}
 Let me make a remark about the $\pm$ or $\mp$ signs. If $\mc F v(a)$ is nonvanishing on $\mathbb R^{-}$, then $b_{\pm n}$ should be read as $b_{+n}$; if $\mc Fv(a)$ is nonvanishing on $\mathbb R^{+}$, then $b_{\pm n}$ should be read as $b_{-n}$. The proof should be read in the same way. \\
 \\
 Proof: Fix  $f(an_t)=\langle \pi_{u, \pm}(an_t) \tau, v \rangle$ bounded on $P_0$ by $C_f$. Suppose that $\mc F v(a)$ is nonvanishing on $\mathbb R^{-}$ or $\mathbb R^{+}$.
 \begin{enumerate}
 \item Suppose that $|f(a n_t)|^2 \leq C_{\mu,f} a^{\mu }$ for  $\mu>0$.  For $-\mu < \epsilon< 0$, the left hand side of Equation (\ref{p0es}) converges. Since $\mc F v(a)$ is nonvanishing on $\mathbb R^{\mp}$, 
 $\int_{0}^{\infty} a^{-\frac{\epsilon}{2}+u_0} | \mc Fv( \mp a)|^2 d a >0$.  Then the sum $\sum_{\frac{1}{\mathfrak p} \leq n } n^
{\frac{\epsilon}{2}-1-u_0} |b_{\pm n}|^2$ becomes a factor and must remain bounded by a constant depending on $f$ and $\epsilon$.
\item Let  $\epsilon>0$, $\delta >0$ and $a_1^2  > \frac{1}{\delta {\mathfrak p}}$. By Prop. \ref{pestimate} we have 
\begin{equation}
\begin{split}
&  (\sum_{\frac{1}{\mathfrak p} \leq n \leq  a_1^2 \delta, n \in {\mathfrak p}^{-1} \mathbb Z} n^{\frac{\epsilon}{2}-1-u_0} |b_{\pm n}|^2 ) \int_{\delta}^{\infty} a^{-\frac{\epsilon}{2}+u_0} |\mc F v(\mp a)|^2  d a \\
\leq &  \int_{\delta}^{\infty} a^{-\frac{\epsilon}{2}+u_0} |\mc F v(\mp a)|^2 (\sum_{\frac{1}{\mathfrak p} \leq n \leq  a_1^2 a, n \in {\mathfrak p}^{-1} \mathbb Z } n^{\frac{\epsilon}{2}-1-u_0} |b_{\pm n}|^2 ) d a \\
\leq &  \int_{\frac{1}{a_1^2 \mathfrak p}}^{\infty} a^{-\frac{\epsilon}{2}+u_0} (\sum_{\frac{1}{\mathfrak p} \leq n \leq  a a_1^2, n \in {\mathfrak p}^{-1} \mathbb Z} n^{\frac{\epsilon}{2}-1-u_0} |b_{\pm n}|^2) |\mc F v (\mp a)|^2  d a \\
 \leq & 2 {\mathfrak p}^{-1} \int_{0}^{a_1} \int_0^{\mathfrak p} |f(a n_t)|^2 d t a^{\epsilon} \frac{ d a}{a} \\
\leq & 2  C_{f} \frac{a_1^{\epsilon}}{\epsilon}
\end{split}
\end{equation}
Now fix a $\delta>0$ such that  $\int_{\delta}^{\infty} a^{-\frac{\epsilon}{2}+u_0} |\mc F v(\mp a)|^2  d a$ is  positive. 
It follows that there exists $C_{\epsilon, f}>0$ such that for any $a_1^2=\frac{k}{ \delta}$,
 $$ \sum_{\frac{1}{\mathfrak p} \leq n \leq k ,n \in {\mathfrak p}^{-1} \mathbb Z } n^{\frac{\epsilon}{2}-1-u_0} |b_{ \pm n}|^2 < 2 C_f^{\prime} \frac{a_1^{\epsilon}}{\epsilon } =2 C_f^{\prime} k^{\frac{\epsilon}{2}} \delta^{-\frac{\epsilon}{2}} \epsilon^{-1} =C_{\epsilon,f, \delta} k^{\frac{\epsilon}{2}}.$$
 Notice that $\delta$ depends on $v$, therefore also on $f$. We can write $c_{\epsilon, f, \delta}$ as $c_{\epsilon, f}$.
 $\Box$
\end{enumerate}
 If $\tau$ is a cuspidal automorphic distribution in a unitary principal series or complementary series representation, then all automorphic forms $f(g)$ will be bounded and rapidly decaying near the cusp at zero. In this situation, the estimates in Theorem \ref{fcoef} were well-known ( \cite{sc}  \cite{br}). The first estimate can also be obtained by observing that the Rankin-Selberg $L(f \times f, s)$ has a pole at $s=1$ for suitable $f$ and the coefficients of the Dirichlet series are all nonnegative (\cite{go}).  If  the (cuspidal) automorphic representation is a discrete series representation, the automorphic distribution $\tau$ can be embedded in $\mc P(u, \pm)^{-\infty}$ for a suitable $u$ and will have its Fourier coefficients supported on ${\mathfrak p}^{-1} \mathbb N$ or $-{\mathfrak p}^{-1} \mathbb N$. Our estimates of Fourier coefficients also follow similarly upon applying the intertwining operator. The details of how to treat the discrete series representations can be found in \cite{sc} \cite{ms1}.   
 
\subsection{ $L^2$-norms of Bounded Periodical Matrix coefficients}
By considering the converse of Theorem \ref{fcoef}, the equations in Prop. \ref{pestimate} also imply the following.
\begin{thm}\label{l2p} Under the same assumption as Proposition \ref{pestimate}, we have the following estimates.
\begin{enumerate}
\item If  $\epsilon< 0$ and
 $\sum_{n \neq 0, n \in {\mathfrak p}^{-1} \mathbb Z } |n|^{\frac{\epsilon}{2}-1-u_0} |b_{n}|^2 <\infty ,$ then there exists positive constant $C_{\epsilon, \tau}$ such that
 $$\int_0^{a_1} \int_0^{\mathfrak p} |f(a n_t)|^2 d t a^{\epsilon} \frac{ d a}{ a} \leq C_{\epsilon, \tau} \sum_{\pm} \int_{\frac{1}{a_1^2 {\f p}}}^{\infty} a^{-\frac{\epsilon}{2}+u_0} \| \mathcal F v(\pm a)|^2 d a.$$
 In particular,
 $$\int_0^{\infty} \int_0^{\mathfrak p} |f(a n_t)|^2 d t a^{\epsilon} \frac{ d a}{ a} \leq C_{\epsilon, \tau} \sum_{\pm} \int_{0}^{\infty} a^{-\frac{\epsilon}{2}+u_0} \| \mathcal F v(\pm a)|^2 d a.$$
 \item If $\epsilon >0$ and $\sum_{|n| \leq k, n \in {\mathfrak p}^{-1} \mathbb Z} |n|^{\frac{\epsilon}{2}-1-u_0} |b_{  n}|^2 < C_{\epsilon,\tau} k^{\frac{\epsilon}{2}}$ for any $k>1$, then 
 $$\int_0^{a_1} \int_0^{\mathfrak p} |f(a n_t)|^2 d t a^{\epsilon} \frac{ d a}{ a} \leq C_{\epsilon, \tau} a_1^{\epsilon} { \mathfrak p} \sum_{\pm}  \int_{\frac{1}{a_1^2 {\f p}}}^{\infty} a^{u_0} \| \mathcal F v(\pm a)|^2 d a.$$
 \commentout{\item If $\epsilon=0$ and $\sum_{|n| \leq k} |n|^{-1-u_0} |b_{ n}|^2 < C_{0,\tau} \ln k (k>1)$,
 then  $$\int_0^{a_1} \int_0^{\mathfrak p} |f(a n_t)|^2 d t \frac{ d a}{ a} \leq C_{\tau} \sum_{\pm} \int_{\frac{1}{a_1^2}}^{\infty} a^{u_0} \| \mathcal F v(\pm a)|^2 (\log a a_1^2) d a.$$}
\end{enumerate}
\end{thm}
We shall remark that this theorem holds even $\mc P(u, \pm)$ is not unitary. \\
\\
Combining Theorems \ref{fcoef} and \ref{l2p}, we have
\begin{cor}[$\epsilon > 0$]\label{dayuning} Under the same assumption as Prop. \ref{pestimate}, suppose for some $\phi \in \mc P(-u, \pm)^{\infty}$ the function $f(an_t)=\langle \pi_{u, \pm}(a n_t) \tau, \phi \rangle$ is bounded on $P_0$ and $\mc F\phi(a)$ is nonvanishing on  both $\mathbb R^+ $ and $\mathbb R^-$. Then for any $\epsilon >0$ and $v \in \mc P(-u, \pm)^{\infty}$, we have
\begin{equation}\label{dayuning1}
\int_0^{a_1} \int_0^{\mathfrak p} |\langle \pi_{u, \pm}(a n_t) \tau, v \rangle |^2 d t a^{\epsilon} \frac{ d a}{ a} \leq C_{\epsilon, \tau} a_1^{\epsilon}  \int_{|a| \geq \frac{1}{a_1^2 \mathfrak p}} |a|^{u_0} \| \mathcal F v( a)|^2 d a .
\end{equation}
In particular, if $\mc P(u, \pm)$ is unitary, we have
$$\int_0^{a_1} \int_0^{\mathfrak p} |\langle \pi_{u, \pm}(a n_t) \tau, v \rangle|^2 d t a^{\epsilon} \frac{ d a}{ a} \leq C_{\epsilon, \tau} a_1^{\epsilon} \| v \|^2_{\mc P(-u, \pm)},$$
$$\int_K \int_0^{a_1} \int_0^{\mathfrak p} |\langle \pi_{u, \pm}(k a n_t) \tau, v \rangle|^2 d t a^{\epsilon} \frac{ d a}{ a} d k \leq C_{\epsilon, \tau} a_1^{\epsilon}  \| v \|^2_{\mc P(-u, \pm)} $$
for every $v \in \mc P(-u, \pm)$.
 \end{cor}
Proof: We only need to prove the second statement. If $u_0=0$, i.e., $\mc P(u, \pm)$ is a unitary principal series, then 
$$\int_{|a| \geq \frac{1}{a_1^2}} |a|^{u_0} \| \mathcal F v( a)|^2 d a \leq \| \mc F v(x) \|^2_{L^2(\mathbb R)} = \| v \|^2_{\mc P(u, \pm)}.$$
If $\mc P(-u, +)$ is a complementary series representation, then the unitary Hilbert norm $\| v \|_{\mc P(-u, \pm)}$ is given by exactly  the square root of
$$\int |x|^{u} \| \mathcal F v( x)|^2 d x ,$$
up to a normalizing factor depending on $u$. Hence we have
$$\int_0^{a_1} \int_0^{\mathfrak p} |f(a n_t)|^2 d t a^{\epsilon} \frac{ d a}{ a} \leq C_{\epsilon, \tau} a_1^{\epsilon} \| v \|^2_{\mc P(-u, \pm)}.$$
Observe that
$$ \langle \pi_{u, \pm}(k a n_t) \tau, v \rangle= \langle \pi_{u, \pm}(a n_t) \tau, \pi_{-u, \pm}(k^{-1} ) v \rangle$$
and $\|\pi(-u, \pm)(k^{-1} ) v \|_{\mc P(-u,\pm)}= \| v \|_{\mc P(-u, \pm)}$. The inequalities in the second statement hold for $v \in \mc P(-u, \pm)^{\infty}$. Therefore, they must  also hold for $v \in \mc P(-u, \pm)$.  $\Box$.\\
\\
Notice that Inequality (\ref{dayuning1}) is true for all $\Re(u)<1$, in particular for $u$ with $\mc P(-u, \pm)$ reducible. Hence it applies to discrete series representation $D_n$. In addition, the norm on the right hand side of Inequality (\ref{dayuning1}) is bounded by 
$$C_{\epsilon, \tau} a_1^{\epsilon}  \int_{a \in \mathbb R} |a|^{u_0} \| \mathcal F v( a)|^2 d a$$
By the Kirillov model, this integral is a constant multiple of the unitary norm $\| v \|_{D_n}$ (\cite{js}).
We have 
\begin{cor}[discrete series case]\label{discrete}
Let $D_n$ be a discrete series representation. Let $\tau$ be a periodic distribution in $D_n^{-\infty}$ with period $\f p$. Suppose that for some $\phi \in D_{-n}^{\infty}$, the function $\langle D_n(a n_t) \tau, \phi \rangle$ is bounded on $P_0$. Then for any $\epsilon>0$ and $v \in D_{-n}^{\infty}$,
$$\int_0^{a_1} \int_0^{\mathfrak p} |\langle D_n(a n_t) \tau, v \rangle|^2 d t a^{\epsilon} \frac{ d a}{ a} \leq C_{\epsilon, \tau} a_1^{\epsilon} \| v \|^2_{D_{-n}},$$
$$\int_K \int_0^{a_1} \int_0^{\mathfrak p} |\langle D_n (k a n_t) \tau, v \rangle|^2 d t a^{\epsilon} \frac{ d a}{ a} d k \leq C_{\epsilon, \tau} a_1^{\epsilon}  \| v \|^2_{D_{-n}} $$
for every $v \in \mc D_{-n}^{\infty}$ and therefore $v \in \mc D_{-n}$. Here $D_{-n}$ is the dual of $D_n$.
\end{cor}
\vspace{0.2 in}
\noindent
Notice that Theorem \ref{l2p} holds for each $\pi_{-u, \pm}(k) v$. We obtain
 \begin{cor}[$\epsilon < 0$]\label{xiaoyuning} Let $\mathcal P(u, \pm)$ be a unitary  representation. Under the assumptions of Prop. \ref{pestimate},
 suppose that $\epsilon< 0$ and
 $\sum_{n \neq 0 } |n|^{\frac{\epsilon}{2}-1-u_0} |b_{ n}|^2 <\infty.$ Then there exists $C_{\epsilon, \tau} >0$ such that
 $$\int_K \int_0^{a_1} \int_0^{\mathfrak p} |\langle \pi_{u, \pm}(k a n_t) \tau, v \rangle |^2 d t a^{\epsilon} \frac{ d a}{ a} d k \leq C_{\epsilon, \tau}  \int_{|x| > \frac{1}{a_1^2 \mathfrak p}} |x|^{-\frac{\epsilon}{2}+u_0} | \mathcal F (\pi_{-u, \pm}(k) v)( x)|^2 d x.$$
 In particular,
 $$\int_K \int \int_0^{\mathfrak p} |\langle \pi_{u, \pm}(k a n_t) \tau, v \rangle |^2 d t a^{\epsilon} \frac{ d a}{ a} d k \leq C_{\epsilon, \tau}  \int_K \int_{-\infty}^{\infty} |x|^{-\frac{\epsilon}{2}+u_0} | \mathcal F (\pi_{-u, \pm} (k) v)( x)|^2 d x d k;$$
 Both inequalities hold for those $v \in \mc P(-u, \pm)$ with which the right hand sides converge.

\end{cor}

In the case of automorphic forms, our $L^2$ norms are estimated over  a Siegel subset, but with the measure $a^{\epsilon} \frac{ d a}{a} d k d t$, while the Siegel set is often equipped with the measure $a^2 \frac{ d a}{a} d k d t$. The bounds we have are certain norms on the representation. This allows us to treat everything at the representation level. If $\epsilon >0$, the bounds come from the Hilbert norm of the automorphic  representation. We have nothing to improve on.  If $\epsilon <0$ , we will need to further study the norm 
\begin{equation}\label{triplenorm}
\vertiii{v}_{\frac{\epsilon}{2}-u_0}=\int_K \int |x|^{-\frac{\epsilon}{2}+ u_0} | \mathcal F (\pi_{-u, \pm} (k) v)( x)|^2 d x d k
\end{equation}
in more details. Our goal is to bound $\vertiii{v}_{\frac{\epsilon}{2}-u_0}$ by a more tangible norm. A natural choice is a norm  coming from the complementary series construction.

\section{$K$-invariant Norms and complementary series}

Let $\Re(u) >-1$. Recall that the smooth vectors in the noncompact picture of unitarizable $\mc P(u, \pm)$ are bounded smooth functions on $\mathbb R$ with integrable Fourier transform. The Fourier transforms  are indeed fast decaying at $\infty$, but singular at zero.
For any bounded smooth function $\phi$ with locally square integrable Fourier transform, let us define 
$$\| \phi \|_{C_u}^2= \int |x|^{-u} | \mathcal F ( \phi )( x)|^2 d x, \qquad (\forall \,\, u \in (-1, 1))$$
whenever such an integral converges. This norm is indeed the unitary norm of the complementary series $C_u$,  upto a normalizing factor. The standard norm $\| * \|_u$ for the complementary series is often constructed using the standard intertwining operator $A_u$ (\cite{knapp}). Our norm $\| * \|_{C_u}$ differs from the $\|* \|_u$ by a normalizing factor. The standard norm $\| * \|_u$ has a pole at $u=0$. The norm $\| * \|_{C_u}$ does not. Hence $\| * \|_{C_u}$ is potentially easier to use.  In this section, we will first review the basic theory of complementary series. Then we will  use $\| \cdot \|_{C_u}$ to bound the norm $\vertiii{ \cdot}_{u}$. Our main references are \cite{knapp} \cite{cas}.
\subsection{Intertwining operator and complementary series}
The standard intertwining operator $A_u: P(u, +)^{\infty} \rightarrow P(-u, +)^{\infty}$ is well-defined for $\Re{u} >0$ and has meromorphic continuation on $\mathbb C$. In the noncompact picture,
$$A_u(f)(x)=\int \frac{f(y)}{|x-y|^{1-u}} d y.$$
Let $\langle *,* \rangle$ be the complex linear $G$-invariant pairing
$$\mc P(u,+) \times \mc P(-u, +) \rightarrow \mathbb C$$
defined by
$$\langle f_1, f_2 \rangle=\int f_1(x) f_2(x) d x \qquad (f_1 \in \mc P(u,+), f_2 \in \mc P(u,-)).$$
For any $\phi, \psi \in \mc P(u, +)$, we define
$$\langle \phi, \psi \rangle_{u}= \langle A_u(\phi), \psi \rangle.$$
This is a $G$-invariant bilinear form on $\mc P(u,+)^{\infty}$. When $u$ is real and $0< u< 1$,
$$( \phi, \psi)_u= \langle A_u(\phi), \overline{\psi} \rangle_u$$
yields an $G$-invariant inner product on $\mc P(u, +)^{\infty}$. Its completion is often called a complementary series representation of $G$, which is  irreducible and unitary. \\
\\ 
In the noncompact picture, the standard basis for the $K$-types of $\mc P(u, +)$ is given by
$$v_{2m}^{(u)}=(1+x^2)^{-\frac{1+u}{2}}(\frac{1+x i}{1-xi})^m \qquad (m \in \mathbb Z).$$
The intertwining operator $A_u$ maps $v_{2m}^{(u)}$ to $c_{2m}^{(u)} v_{2m}^{-u}$. The constant 
$$c_{2m}^{(u)}=\frac{ (-1)^m 2^{1-u} \pi \Gamma(u)}{\Gamma(\frac{u+1}{2}+m)\Gamma(\frac{u+1}{2}-m)}=\frac{2^{1-u}  \Gamma(u) \Gamma(m+\frac{-u+1}{2}) \sin(\frac{u+1}{2}\pi)}{\Gamma(\frac{u+1}{2}+m) } .$$
See \cite{cas}. We make two observations here. First, the formula above in fact uniquely determined the analytic continuation of the intertwining operator $A_u$. Secondly, for  $u \notin 2 \mathbb Z+1$,
$$\frac{\Gamma(\frac{-u+1}{2}+m)}{\Gamma(\frac{u+1}{2}+m)} \sim c_u m^{-u} \qquad (m \rightarrow \infty).$$ We have
\begin{lem}\label{com-norm}
 For a fixed $u \in (-1,0)$ or $u \in (0,1)$, there exist  positive constants $c_{u}$, $c_u^{\prime}$ such that
$$c_u^{\prime} (1+|m|)^{-u} \leq (v_m^{(u)}, v_m^{(u)})_u \leq c_u (1+|m|)^{-u} \qquad ( m \in \mathbb Z)).$$
\end{lem}
The intertwining operator $A_u$ has a pole at $u=0$. Hence we must exclude $u=0$ from our estimates.

\subsection{Normalizing  $(*,*)_u$}
Recall that for $u \in (0,1)$
$$(\phi, \psi)_u=\int \int \frac{\phi(x) \overline{\psi(y)}}{|x-y|^{1-u}} d x d y \qquad (\phi, \psi \in \mc P(u, +)^{\infty}),$$
and 
$$(\phi, \psi)_{C_u}=\int |\xi|^{-u} \mc F(\phi)(\xi) \overline{\mc F(\psi)(\xi)} d \xi.$$
By Fourier inversion formula, we have
$$(\phi, \psi)_{C_u}= G(u) (\phi, \psi)_u,$$
where $\int |\xi|^{-u} \exp -2\pi i x \xi d \xi= G(u) |x|^{-1+u}$ . This is true for $u \in (0,1)$ and can be analytically continued to $u \in (-1,1)$, since the function $G(u)$ can be expressed in terms of $\Gamma$-functions and possesses a zero at $u=0$ (\cite{ms}). Hence we have
 $$\| v_m^{(u)} \|_{C_u}^2=G(u) \| v_m^{(u)}\|_u^2$$ for $u \in (-1, 1)$. By Lemma \ref{com-norm} we have the following estimates:
\begin{thm}\label{com-norm1} For $u \in (-1,0]$, there exist positive constants $q_u$, $q_u^{\prime}$ depending continuously on $u$ such that $q_0=q_0^{\prime}=1$ and
$$q_u^{\prime} (1+m^2)^{-\frac{u}{2}} \leq \| v_{2m}^{(u)}\|_{C_u}^2 \leq q_u (1+m^2)^{-\frac{u}{2}}.$$
\end{thm}

\subsection{Bounds by the complementary norm: $\mathcal P(i \lambda, +)$ case}
Fix $v \in \mc P(i \lambda, +)^{\infty}$ with $\lambda \in \mathbb R$.
Recall that we are interested in the norm
$$ \vertiii{ v}_{u}= \int_K \int |x|^{-u} | \mathcal F (\pi_{ i \lambda, +} (k) v)( x)|^2 d x d k=\int_K \| \pi_{i \lambda, +}(k) v \|_{C_{u}}^2 d k \qquad (u <0).$$
Clearly, this norm is $K$-invariant. Hence we will need to estimate
$\vertiii{ v_{2m}^{(i \lambda)}}_{u}= \| v_{2m}^{(i \lambda)} \|_{C_u}$.
\begin{thm}\label{ktype} Let $ u \in (-1,0)$. Then there exists a positive constant $c_{u}$ such that $\forall \, m \in \mathbb Z$
$$\vertiii{ v_{2m}^{(i \lambda)}}_{u}^2 \leq c_{u} (1+ |m|^{-u}).$$
\end{thm}
Proof: Observe that
$$v_{2m}^{(i \lambda)}(x)=(1+x^2)^{\frac{-i\lambda+u}{2}}v_{2m}^{(u)}.$$
Under the compact picture of $\mathcal P(u, +)$, $v_{2m}^{(i \lambda)}$ becomes
$$|\sin \theta|^{i \lambda-u} \exp 2m i \theta, \qquad (\cot \theta=x).$$
The function $|\sin \theta|^{i \lambda-u}$ has period $\pi$ and $L^1$ derivative. Hence its Fourier series expansion
$$\sum_{k \in \mathbb Z} a_{2k} \exp 2k i \theta,$$
satisfy that $|a_{2k}| \leq h_{u} (1+ k^2)^{-\frac{1}{2}}$ for some positive constant $h_u$.
We obtain
$$v_{2m}^{(i \lambda)}=\sum_{k \in \mathbb Z} a_{2k} v_{2m+2k}^{(u)}.$$
It follows that
$$\| v_{2m}^{(i \lambda)} \|_{C_u}^2= \sum_{k \in \mathbb Z} |a_{2k}|^2 \|v_{2m+2k}^{(u)} \|_{C_u}^2 \leq h_{u}^2 q_{u} \sum_{k \in \mathbb Z} \frac{(1+(m+k)^2)^{-\frac{u}{2}}}{k^2+1}  \leq h_{u}^2 q_{u} \sum_{k \in \mathbb Z} \frac{(1+2m^2)^{-\frac{u}{2}}(1+2k^2)^{-\frac{u}{2}}}{k^2+1} ,$$
which will be bounded by a multiple of  $(1+m^2)^{-\frac{u}{2}}$.
$\Box$ \\
\\
For $u \in (-1, 0)$ the map
$$ v(x) \in \mc P(i \lambda,+)^{\infty} \rightarrow (1+x^2)^{\frac{ i \lambda-u}{2}} v(x) \in \mc P(u, +)^{\infty}$$
preserves the $K$ action and maps $v_{2m}^{(i \lambda)}$ to $v_{2m}^{(u)}$. By Theorem \ref{com-norm1} and \ref{ktype}
there is a constant $c_{u}$ such that
$$\vertiii{v_{2m}^{(i \lambda)}}_u \leq c_{u} \| v_{2m}^{(u)} \|_{C_u}  $$
We have
\begin{thm}\label{exotic} For $u \in (-1,0)$ and $\lambda \in \mathbb R$, there exists a positive constant $c_{u}$ such that
$$\vertiii{ v(x)}_{u}^2 \leq c_{u} \| (1+x^2)^{\frac{ i \lambda-u}{2}}  v(x) \|_{C_u} \qquad (\forall \,\,\, v(x) \in \mc P(i \lambda, +)^{\infty}).$$
\end{thm}
\vspace{0.1in}
\noindent
Under the assumption of  Cor \ref{xiaoyuning}, applying Theorem \ref{exotic}, $\langle \pi_{u, \pm}(k a n_t) \tau, v \rangle$ will be in $L^2(G/N_{\mathfrak p}, a^{\epsilon} \frac{d a}{a} d t d k)$ as long as $\| (1+x^2)^{\frac{ 2i \lambda-\epsilon}{4}}  v(x) \|_{C_{\frac{\epsilon}{2}}}$ is bounded with $\epsilon \in (-2,0)$.
\subsection{Bounds by the complementary norm: $\mathcal P(i \lambda, -)$ case}
Let $u \in (-1, 0)$ and $\lambda \in \mathbb R$. The $K$-types in $\mc P(i \lambda, -)$ are
$$v_{2m+1}^{(i \lambda)}(x)=(1+x^2)^{\frac{-i \lambda+u}{2}} (\frac{1+xi}{1-xi})^{m+\frac{1}{2}} \qquad (m \in \mathbb Z).$$
Here $x= \cot \theta \,\,\,(\theta \in (0, \pi))$ and $(\frac{1+xi}{1-xi})^{m+\frac{1}{2}}= \exp i (2m+1) \theta$ is well-defined. Our goal is to estimate $\|v_{2m+1}^{(i \lambda)}(x)\|_{C_u}$. We still have
$$v_{2m+1}^{(i \lambda)}(x)=v_{2m}^{(u)} (\frac{1+xi}{1-xi})^{\frac{1}{2}} (1+x^2)^{\frac{u- i \lambda}{2}}.$$
In the compact picture of $\mc P(u,+)$, $v_{2m+1}^{(i \lambda)}(x)$ becomes $\sgn (\sin \theta) |\sin \theta|^{-u+i \lambda} \exp (2 m+1)i \theta$.
Notice that this function has period $\pi$ and take the same value as $|\sin \theta|^{-u+i \lambda} \exp (2 m+1)i \theta$ when $\theta \in [0, \pi]$.
Observe that $\sgn (\sin (\theta+\pi)) |\sin (\theta+\pi))|^{-u+i \lambda}= -\sgn (\sin \theta) |\sin \theta|^{-u+i \lambda}$. Let $\sgn (\sin \theta) |\sin \theta|^{-u+i \lambda}=\sum_{k \in \mathbb Z} b_{2k-1} \exp (2k-1) i \theta$ be its Fourier expansion. Again, the function $\sgn (\sin \theta) |\sin \theta|^{-u+i \lambda}$ has $L^1$-derivative. Hence $|b_{2k-1}| \leq c_{u} \frac{1}{|2k-1|}$. By a similar argument as $\mc P(i \lambda, +)$ case we have

\begin{thm}\label{ktype1} Let $ u \in (-1,0)$. Then there exists a positive constant $c_{u}$ such that
$$  \vertiii{ v_{2m+1}^{(i \lambda)}}_{u}^2 \leq c_{u} (1+ |m|^{-u}).$$
\end{thm}
\begin{thm}\label{exotic1} For $u \in (-1,0)$ and $\lambda \in \mathbb R$, there exists a positive constant $c_{u}$ such that
$$\vertiii{ v(x)}_{u}^2 \leq c_{u} \| (1+x^2)^{\frac{ i \lambda-u}{2}} (\frac{1+xi}{1-xi})^{\frac{1}{2}} v(x) \|_{C_u} \qquad (v(x) \in \mc P(i \lambda, -)^{\infty}).$$
\end{thm}
Proof: Consider the map
$$I: \mc P(i \lambda, -)^{\infty} \rightarrow \mc P(u, +)^{-\infty}$$
defined by
$$I(v)(x)=(1+x^2)^{\frac{ i \lambda-u}{2}} (\frac{1+xi}{1-xi})^{\frac{1}{2}}v(x).$$
$I$ maps the orthogonal basis $\{ v_{2m-1}^{(i \lambda)}: m \in \mathbb Z \}$ of $\vertiii{*}$ to orthogonal basis $\{ v_{2m}^{(u)}: m \in \mathbb Z \}$ 
of the complementary series $C_u$. In addition, one can easily check that $I$ is bounded. Our theorem then follows. Contrary to the spherical case, the operator $I$ is no longer $K$-invariant. $\Box$
\subsection{Bounds by the complementary norm: $\mc P(u, +)$ case}
Let $u \in (-1,1)$. Then $\mc P(u, +)$ is the complementary series $C_{u}$. For $\mu  <0$ and $ v \in \mc P(u, +)^{\infty}$, we are interested in
$$\vertiii{ v(x)}_{u+\mu}^2 = \int_K \int_{\mathbb R} |x|^{-u-\mu} |\mc F(\pi_{u, +}(k)v)(x)|^2 d x d k.$$
For our purpose, we will assume that $u + \mu> -1$.

\begin{thm}
Let $u \in (-1, 1)$ and $\mu \in (-1-u, 0)$. Then  there exists a positive constant $c_{\mu, u}$ such that
$$\vertiii{ v_{2m}^{(u)}}_{u+\mu}^2 \leq c_{\mu, u} (1+ |m|)^{-u-\mu}.$$
\end{thm}

If $u+ \mu \leq 0$, our proof is similar to the proof of Theorem \ref{ktype}. If $0 < u+\mu <1$, the proof will be different. We will be a little sketchy.\\
\\
Proof: We have $v_{2m}^{(u)}= v_{2m}^{(u+\mu)}(1+x^2)^{\frac{\mu}{2}}$. Under the compact picture, $v_{2m}^{(u)}= v_{2m}^{(u+\mu)}| \sin \theta|^{-\frac{\mu}{2}}$. Let $ \sum_{k \in \mathbb Z} a_{2k} \exp 2k i \theta$ be the Fourier expansion of $|\sin \theta|^{-\frac{\mu}{2}}$. Since 
$|\sin \theta|^{-\frac{\mu}{2}}$ has $L^1$-derivative, we must have
$|a_{2k}| \leq h_{\mu} (1+ k^2)^{-\frac{1}{2}}$ for a positive constant $h_{\mu}$.
We obtain
$$v_{2m}^{(u)}=\sum_{k \in \mathbb Z} a_{2k} v_{2m+2k}^{(u+ \mu)}.$$
Notice $u+\mu >-1$. If $u+ \mu \leq 0$, by Theorem \ref{com-norm1},
$$\vertiii{ v_{2m}^{(u)}}_{u+\mu}^2=\| v_{2m}^{(u)} \|_{C_{u+ \mu}}^2= \sum_{k \in \mathbb Z} |a_{2k}|^2 \|v_{2m+2k}^{(u+\mu)} \|_{C_{u+\mu}}^2 \leq h_{\mu}^2 q_{u+\mu} \sum_{k \in \mathbb Z} \frac{(1+(m+k)^2)^{-\frac{u+\mu}{2}}}{k^2+1}  $$
$$\leq h_{\mu}^2 q_{u+\mu} \sum_{k \in \mathbb Z} \frac{(1+2m^2)^{-\frac{u+\mu}{2}}(1+2k^2)^{-\frac{u+\mu}{2}}}{k^2+1} ,$$
which will be bounded by a multiple of  $(1+m^2)^{-\frac{u+\mu}{2}}$.\\
\\
If $u+\mu >0$ and $m \neq 0$, we have
$$\sum_{k \in \mathbb Z} \frac{(1+(m+k)^2)^{-\frac{u+\mu}{2}}}{k^2+1}=\sum_{|k| > \frac{|m|}{2}} \frac{(1+(m+k)^2)^{-\frac{u+\mu}{2}}}{k^2+1}+ \sum_{|k| \leq \frac{|m|}{2}} \frac{(1+(m+k)^2)^{-\frac{u+\mu}{2}}}{k^2+1}.$$
The first sum is bounded by $\sum_{|k| > \frac{|m|}{2}} \frac{1}{k^2+1} \leq c |m|^{-1} \leq c |m|^{-u-\mu}$, since $u+ \mu <1$. The second sum is bounded by $c^{\prime} |m|^{-u-\mu}$. We see that
$\| v_{2m}^{(u)} \|_{C_{u+ \mu}}^2 \leq c_{u, \mu} (1+ |m|)^{-u-\mu}.$ $\Box$\\
\\
\noindent
By essentially the same proof as Theorem \ref{exotic}, we have
\begin{thm}\label{exotic2} For $u \in (-1,1)$ and $\mu \in (-1-u, 0)$, there exists a positive constant $c_{u, \mu}$ such that
$$\vertiii{ v(x)}_{u+\mu}^2 \leq c_{u, \mu} \| (1+x^2)^{\frac{ -\mu}{2}}  v(x) \|_{C_{u+\mu}} \qquad (v(x) \in \mc P(u, +)^{\infty}).$$
\end{thm}

\section{$K$-invariant Norms over $G/\Gamma$}
 Let $\Gamma$ be a nonuniform lattice in $SL(2, \mathbb R)$. Then $G/\Gamma$ has a finite volume and a finite number of cusps, $z_1, z_2, \ldots, z_l$. Write
 $G/\Gamma$ as the union of Siegel sets  $S_1, S_2, \ldots S_l$ with a compact set $C_0$ (\cite{borel}). Since $\Gamma$ action is on the right, our standard Siegel set will be near $0$, not $\infty$.  Let $ d g= a\,  d a \,d t \, d k$ be the invariant measure of $G$ under the $KAN$ decomposition. Over each Siegel set $S_i$, the invariant measure can be written as $d g= a_i d a_i d t_i d k$. 
\begin{thm}\label{main1} Let $\Gamma$ be a nonuniform lattice in $SL(2, \mathbb R)$. Let $\mathcal H \subseteq L^2(G/\Gamma)$ be a cuspidal automorphic representation of type $\mc P(-u, \pm)$. Given any $K$-invariant measure $\nu$ on $G/\Gamma$ such that $\nu$ is bounded by $ dg$ on $C_0$ and bounded by $a_i^{\epsilon} \frac{d a_i}{a_i} d t_i d k$ on $S_i$,  there exists a constant $C$ depending on $\nu$ (hence on $\epsilon$) and $\mc H$ such that
\begin{enumerate}
\item If $\epsilon >0$, then 
$$ \| f \|_{L^2(G/\Gamma, d \nu)} \leq C \|f \|_{L^2(G/\Gamma, d g)}, \qquad ( f \in \mc H);$$
\item If $\epsilon <0$, then
\commentout{
$$ \| f \|_{L^2(S_i, d \nu)} \leq C \vertiii{f}_{\frac{\epsilon}{2}-u_0}, \qquad ( f \in \mc H^{\infty}).$$
where $\vertiii{f}_{\frac{\epsilon}{2}-u_0}$ is defined intrinsically for $\mc H^{\infty} \cong \mc P(u, \pm)^{\infty}$ for each cusp $z_i$. Generally, } 
for any $f \in \mc H^{\infty} \cong \mc P(u, \pm)^{\infty}$,
$$\| f \|_{L^2(G/\Gamma, d \nu)} \leq C_{\mc H} \vertiii{f}_{\frac{\epsilon}{2}-u_0}$$  and $\vertiii{f}_{\frac{\epsilon}{2}-u_0}$ will be  bounded  the complementary norm  given in Theorems \ref{exotic} \ref{exotic1} \ref{exotic2}.
\end{enumerate}
\end{thm}
We shall remark that our theorem can be generalized to all nonuniform lattice of a finite covering of $SL(2, \mathbb R)$. \\
\\
Proof:  Let $v \in \mathcal P(-u, \pm)^{\infty}$ and $\sigma \in \mathcal P(u, \pm)^{-\infty}$. Let $f(kan_t)=\langle \pi_{u, \pm}(k a n_t)  \sigma, v \rangle$. Then for any $h \in G$, the left action $$L(h) f(g)=f(h^{-1} g)=\langle \pi_{u, \pm}(h^{-1} g) \sigma, v \rangle= \langle \pi_{u, \pm}( g) \sigma,\pi_{-u, \pm}(h) v \rangle.$$
We see that  the left action on $f(kan_t)$ is equivalent to the action of $\mc P(-u, \pm)$ on $v$.  Fix $\mathcal H \subseteq L^2(G/\Gamma)$,  a cuspidal automorphic representation of type $\mc P(-u, \pm)$. By \cite{sc} \cite{br}, there exists a $\Gamma$-invariant distribution $\tau \in \mathcal P(u, \pm)^{-\infty}$ such that all smooth vectors in $\mc H^{\infty}$ can be written as $\langle \pi_{u, \pm}(g)  \tau, v \rangle$ for some $v \in \mathcal P(-u, \pm)^{\infty}$. \\
\\
Fix $\epsilon > 0$. For each cusp $z_i$, we can use the action of $k_i$ so that $k_i z_i=0$. In the language of Harish-Chandra, this amounts to choose a cuspidal pair $(P, A)$.  By Cor. \ref{dayuning}, for each cusp $z_i$, we can choose a Siegel set $S_i$ and find a constant $C_i$ such that
$$ \| \langle \pi_{u, \pm}(g)  \tau, v \rangle \|_{L^2(S_i, a_i^{\epsilon} \frac{d a_i}{a_i} d t_i d k)} \leq c_i \| v \|_{\mc P(-u, \pm)}=c_i^{\prime} \| \langle \pi_{u, \pm}(g)  \tau, v \rangle \|_{L^2(G/\Gamma)}.$$
Obviously, for the compact set $C_0$, $$\| \langle \pi_{u, \pm}(g)  \tau, v \rangle \|_{L^2(C_0, d g)} \leq \| \langle \pi_{u, \pm}(g)  \tau, v \rangle \|_{L^2(G/\Gamma)}.$$
Hence, our first inequality follows. \\
\\
Fix $\epsilon < 0$. By Cor. \ref{xiaoyuning}, $\| \langle \pi_{u, \pm}(g)  \tau, v \rangle  \|_{L^2(S_i, d \nu)} \leq C \vertiii{v}_{\frac{\epsilon}{2}-u_0}$ defined for each cusp $z_i$. In the cases of $\mc P(-i\lambda, +)$,  By Theorem \ref{exotic}, the norm 
$$\vertiii{v}_{\frac{\epsilon}{2}} \leq C_i \| (1+x^2)^{\frac{2 i \lambda-\epsilon}{4}}  v(x) \|_{C_{\frac{\epsilon}{2}}}.$$
Observe that the map from $\mc P(-i\lambda,+)^{\infty}$ to $\mc P(\frac{\epsilon}{2}, +)^{\infty}$ defined by
$$ v(x) \rightarrow (1+x^2)^{\frac{2 i \lambda-\epsilon}{4}}  v(x)$$
is $K$-invariant and the $\| * \|_{C_{\frac{\epsilon}{2}}}$ is independent of the choices of the unipotent subgroup $N$. Hence 
$\| (1+x^2)^{\frac{2 i \lambda-\epsilon}{4}}  v(x)\|_{C_{\frac{\epsilon}{2}}}$ remains the same for different choices of cusps. Over $C_0$, we have
$$\| \langle \pi_{i \lambda, +}(g)  \tau, v \rangle \|_{L^2(C_0, d g)} \leq \| \langle \pi_{ i\lambda, +}(g)  \tau, v \rangle \|_{L^2(G/\Gamma)}=
 c_2 \| v \|_{\mc P(-i\lambda, +)}  \leq c_2 \| (1+x^2)^{\frac{2 i \lambda-\epsilon}{4}}  v(x)\|_{C_{\frac{\epsilon}{2}}}.$$
We obtain
$$\| \langle \pi_{i \lambda, +}(g)  \tau, v \rangle  \|_{L^2(G/\Gamma, d \nu)} \leq 
C_{\mc H} \|   (1+x^2)^{\frac{2 i \lambda-\epsilon}{4}}  v(x)\|_{C_{\frac{\epsilon}{2}}}.$$
The complementary series case $\mc P(u, +)$ is similar. The nonspherical unitary principal series $\mc P(i \lambda, -)$ is more delicate. Essentially, norms  
$\vertiii{v}_{\frac{\epsilon}{2}}$ with respect to different $N_i$ will be mutually bounded. Hence we still have
$$\| \langle \pi_{i \lambda, +}(g)  \tau, v \rangle  \|_{L^2(G/\Gamma, d \nu)} \leq 
C_{\mc H} \|   {(1+x^2)^{\frac{2 i \lambda-\epsilon}{4}} (\frac{1+xi}{1-xi})^{\frac{1}{2}} v(x)} \|_{C_{\frac{\epsilon}{2}}}.$$
$\Box$
\subsection{Bounds with respect to $\Omega A$}
The $KAN$ decomposition fits naturally in the theory of Fourier-Whittaker coefficients of automorphic forms. It is used by number theorists to conduct analysis on automorphic forms, often over a Siegel set. However to understand the $L$-function of automorphic representation, in particular, the growth of L-function, the natural choice seems to be the $KNA$ decomposition. Both $KAN$ and $KNA$ originated in the Iwasawa decomposition and are closely related to Cartan decomposition. The analysis based on these decomposition seems to be of different flavor and have different implications.  The $G$-invariant measure with respect to $KAN$ decomposition is $a^2 \frac{d a}{a} d n d k$ or $a^{-2} \frac{ d a}{a} d a dn dk$ depending on the choices of $N$. The $G$-invariant measure with respect to $KNA$ decomposition is simply $d k \, d n \, d a$. \\
\\
 Recall that  L-function for a cuspidal automorphic representation of $SL(2, \mathbb R)$ can be represented by a zeta integral over $MA \cong GL(1)$. Hence it is desirable to have an estimate of the $L^2$-norm of automorphic forms over $\Omega A$, where $\Omega$ a compact set with finite measure in $KN$. 

\begin{thm}\label{main2} Let $\Gamma$ be a nonuniform lattice  in $SL(2)$. Suppose that $w \in \Gamma$ and $N_{\f p} \subseteq \Gamma$. Let $\mc H$ be a cuspidal  automorphic representation of $G$ of type $\mc P(i \lambda, \pm)$.  Then there exists a positive constant $C$ depending on $\epsilon, \mc H$ and $T_1$ such that
$$\| f \|_{T_1, \epsilon} \leq C \vertiii{f}_{-\frac{|\epsilon|}{2}} \qquad (f \in \mc H^{\infty}).$$
\end{thm}

Proof: By Theorem \ref{main}, $$\| f \|_{T_1, \epsilon}^2  \leq c_{T_1,\epsilon, \mathfrak p} \int_K \int_0^{\sqrt{1+T_1^2}} ( a^{\epsilon}+a^{-\epsilon}) \int_0^{\mathfrak p} | f (ka n_t)|^2 d t \frac{ d a}{ a} d k $$
$$ \leq C_{T_1, \epsilon, \f p} ((1+T_1^2)^{\epsilon}+1) \int_K \int_0^{\sqrt{1+T_1^2}} a^{-|\epsilon|} \int_0^{\mathfrak p} | f (ka n_t)|^2 d t \frac{ d a}{ a} d k.$$
Since $\mc H$ is cuspidal, the $K$-finite functions in $\mc H$ are bounded and rapidly decaying near the cusp $0$. Again, we write 
$f(g) \in \mc H$ as matrix coefficient $\langle \pi_{i \lambda, \pm}(g)  \tau, v \rangle$ for some $ v \in \mc P(-i\lambda, \pm)$ and $\tau \in \mc P(\i \lambda, \pm)^{-\infty}$. Obviously, $\tau$ will have no constant term in Fourier expansion. Its Fourier coefficients have the convergence specified in Theorem \ref{fcoef}.
By Cor \ref{dayuning} \ref{xiaoyuning}, there exists $ C_{\epsilon, \mc H, T_1} >0$ such that
$$\| f \|_{T_1, \epsilon}^2 \leq C_{\epsilon, \mc H, T_1}  \vertiii{f}_{-\frac{|\epsilon|}{2}}^2 \qquad (f \in \mc H^{\infty}).$$
Our theorem then follows. $\Box$ 

\begin{cor}\label{main3}   Let $\Gamma$ be a nonuniform lattice  in $SL(2, \mathbb R)$. Suppose that $w \in \Gamma$ and $N_{\f p} \subseteq \Gamma$. Let $\mc H$ be a cuspidal  automorphic representation of $G$ of type $\mc P(i \lambda, \pm)$.   Let $\Omega$ be a compact 2 dimensional domain in $KN$.  Let $\epsilon \in \mathbb R$. Then there exists a positive constant $C$ depending on $\epsilon, \mc H$ and $\Omega$ such that
$$\| f \|_{L^2(\Omega A, a^{\epsilon} \frac{ d \, a }{a} d t d k)} \leq C \vertiii{f}_{-\frac{|\epsilon|}{2}} \qquad (f \in \mc H^{\infty}).$$
\end{cor}
Proof: Obviously, any compact set $\Omega$ in $KN$ is contained in some $K N_{T_1}$. Hence $\Omega A \subseteq X_{T_1}$. Then our assertion follows from the previous theorem. $\Box$ \\
\\
We shall remark that our results also apply to cuspidal automorphic representations of type $\mc P(-u, +)$ with $u \in (-1,1)$. The bound will be a constant multiple of $\vertiii{f}_{-u-\frac{|\epsilon|}{2}}$ as in Theorem \ref{main2}.

\subsection{Applications to Unitary Eisenstein series}
 We shall remark that the Theorem \ref{main2} remains to be true if 
 \begin{enumerate}
 \item $w \in \Gamma$ and $N_{\f p} \subseteq \Gamma$;
 \item the Fourier coefficient $b_n$ of $\tau$ satisfies the conditions that $b_0=0$ and $\sum |n|^{-\frac{\epsilon}{2}-1-u_0} |b_{  n}|^2 <\infty$ for $\epsilon>0$.
 \end{enumerate}
The following proposition follows directly from Theorem \ref{l2p}.
  \begin{prop}
  Let $\Gamma$ be a discrete subgroup of $SL(2)$ such that $w \in \Gamma$ and $N_{\f p} \subseteq \Gamma$. Let $\mc V$ be an automorphic representation of type
  $\mc P(i \lambda, \pm)$. In addition, we can assume $\mc V$ is given by $\langle \pi_{i \lambda, \pm}(g) \tau, v \rangle$ with $\tau \in \mc P(i\lambda, \pm)^{-\infty}$. Let $\epsilon > 0$ and suppose  $\tau=\sum_{n \in {\mathfrak p}^{-1} \mathbb Z, n \neq 0}^* b_{n} \exp  2 \pi i x n$ with  $\sum_{|n| \leq k} |n|^{-\frac{\epsilon}{2}-1} |b_{  n}|^2 <\infty$. Then 
  $$\| \langle \pi_{i \lambda, \pm}(g) \tau, v \rangle \|_{T_1, \epsilon} \leq C \vertiii{v}_{-\frac{\epsilon}{2}} \qquad (v \in \mc H^{\infty}).$$
  \end{prop}
  
 If $\Gamma$ is a congruence subgroup containing $w$ and   the unitary Eisenstein series is cuspidal at $0$ and $\infty$, we have   
\begin{cor}  Let $\Gamma$ be a congruent subgroup of $SL(2, \mathbb R)$ such that $w \in \Gamma$. Let $\mc V$ be an Eisenstein series of type $\mc P(i \lambda, \pm)$ and $\epsilon \in \mathbb R$. Suppose that $\mc V$ has zero constant term with respect to $N$.   Then 
  $$\| f \|_{T_1, \epsilon} \leq C \vertiii{f}_{-\frac{|\epsilon|}{2}} \qquad (f \in \mc V).$$
\end{cor}
Proof:  The Fourier coefficients of Eisenstein series for congruence subgroups are computable (\cite{go}). It can be checked that $\sum |n|^{-\frac{\epsilon}{2}-1} |b_{  n}|^2 <\infty$ for $\epsilon>0$. $\Box$

\end{document}